\documentclass[letterpaper, 10 pt, conference]{ieeeconf}

\IEEEoverridecommandlockouts  

\newtheorem{theorem}{Theorem}[section]
\newtheorem{lemma}[theorem]{Lemma}
\newtheorem{definition}{Definition}
\newtheorem{corollary}[theorem]{Corollary}

\newtheorem{remark}{Remark}
\newtheorem{assumption}{Assumption}
\newtheorem{example}{Example}

\newcommand*{\QEDB}{\hfill\ensuremath{\square}}
\newcommand*{\QEDBL}{\hfill\ensuremath{\blacksquare}}
\newcommand\oprocendsymbol{\hbox{$\square$}}
\newcommand\oprocend{\relax\ifmmode\else\unskip\hfill%
\fi\oprocendsymbol}

\newcommand{\cD}{{\cal D}}

\newcommand{\cS}{{\cal S}}

\newcommand{\cX}{{\cal X}}

\newcommand{\bc}{{\bf c}}
\newcommand{\ba}{{\bf a}}
\newcommand{\bb}{{\bf b}}

\newcommand{\be}{{\bf e}}

\newcommand{\br}{{\bf r}}
\newcommand{\bs}{{\bf s}}
\newcommand{\bx}{{\bf x}}
\newcommand{\bu}{{\bf u}}
\newcommand{\bv}{{\bf v}}
\newcommand{\bw}{{\bf w}}

\newcommand{\bz}{{\bf z}}
\newcommand{\by}{{\bf y}}
\newcommand{\bA}{{\bf A}}
\newcommand{\bB}{{\bf B}}

\newcommand{\bG}{{\bf G}}

\newcommand{\bH}{{\bf H}}

\newcommand{\bI}{{\bf I}}

\newcommand{\bU}{{\bf U}}

\newcommand{\bV}{{\bf V}}
\newcommand{\bbP}{{\mathbb{P}}}
\newcommand{\bbR}{{\mathbb{R}}}

\newcommand{\norm}[1]{\Vert #1 \Vert}

\newcommand{\pr}[1]{\mathbb{P}\left[#1\right]}
\newcommand{\ev}[1]{\mathbb{E}\left[#1\right]}

\newcommand{\sw}{\mathrm{subW}}

\usepackage{hyperref}
\usepackage{graphicx,color}
\graphicspath{{./IMG/}{images/}{./IMG/SAGEMATH_PLOTS/}}
\usepackage{amsmath}
\usepackage{amssymb}
\usepackage{mathrsfs}
\usepackage{subfigure}
\usepackage{url}
\usepackage{booktabs}
\usepackage{array}
\usepackage[table]{xcolor}
\usepackage{mathtools} 		
\usepackage{epstopdf}
\usepackage{cite}
\usepackage[vlined,ruled]{algorithm2e}
\usepackage{upgreek}
\usepackage{dsfont}

\title{\LARGE \bf
Online Stochastic Gradient Methods Under Sub-Weibull Noise and the Polyak-\L ojasiewicz Condition
}

\author{Seunghyun Kim, Liam Madden, and Emiliano Dall'Anese
\thanks{S. Kim and L. Madden are with the Department of Applied Mathematics at the University of Colorado Boulder; E. Dall'Anese is with the Department of Electrical, Computer and Energy Engineering and with the Department of Applied Mathematics at the University of Colorado Boulder. }
\thanks{This work was supported by the National Science Foundation (NSF) CAREER award 1941896 and the NSF award 1923298.}
}

\begin{document}
\maketitle

\allowdisplaybreaks[3]

\begin{abstract}
This paper focuses on the online gradient and proximal-gradient methods with stochastic gradient errors. In particular, we examine the performance of the online gradient descent method when the cost satisfies the Polyak-\L ojasiewicz (PL) inequality. We provide bounds in expectation and in high probability (that hold iteration-wise), with the latter derived by leveraging a sub-Weibull model for the errors affecting the gradient.
The convergence results show that the instantaneous regret converges linearly up to an error that depends on the variability of the problem and the statistics of the sub-Weibull gradient error. Similar convergence results are then provided for the online proximal-gradient method, under the assumption that the composite cost satisfies the proximal-PL condition. In the case of static costs, we provide new bounds for the regret incurred by these methods when the gradient errors are modeled as sub-Weibull random variables. Illustrative  simulations are provided to corroborate the technical findings.
\end{abstract}

\section{Introduction}

This paper considers online gradient descent and the online proximal-gradient  methods for dynamic optimization and learning~\cite{popkov2005gradient, towfic2013distributed,bedi2018tracking,selvaratnam2018Numerical,yi2016tracking, mokhtari2016online,chang2021online,besbes2015non}. Because of their computational tractability, these are attractive first-order methods for solving a number of learning and optimization tasks where  data  points  and functions are  processed  on-the-fly  and  without   storage. Online gradient and proximal-gradient descent are powerful methods also in the context of online stochastic optimization~\cite{shames2020online,cao2020online}, stochastic learning~\cite{vlaski2020tracking,hallak2020regret}, and feedback-based optimization~\cite{hauswirth2021optimization,ospina2022feedback}.

We examine the performance of online gradient and proximal-gradient descent in the presence of \emph{inexact} gradient information, and when the cost to be minimized satisfies the \emph{Polyak-\L ojasiewicz (PL) condition}~\cite{karimi2016linear}. Formally, we consider  a optimization problem of the form
\begin{align}
\label{eq:main-problem}
\min_{\bx \in \mathbb{R}^n} F_t(\bx) := f_t(\bx) + g_t(\bx)
\end{align}  
where $t \in \mathbb{N}$ is the time index, $f_t: \mathcal{D} \rightarrow \mathbb{R}$ is a continuously differentiable function with a Lipschitz-continuous gradient at each time $t$, $\mathcal{D} \subseteq \mathbb{R}^n$ an open and non-empty convex set, and $g_t: \mathcal{D}  \to \mathbb{R} \cup \{ +\infty \}$ is a  closed, convex and proper function for all $t$, possibly not differentiable. Accordingly, we consider two main cases: 

\noindent  \emph{c1)} $g_t(\bx) \equiv 0$, $\bx \mapsto f_t(\bx)$ satisfies the PL inequality for all $t$, and an inexact gradient is available; and, 

\noindent \emph{c2)} $\bx \mapsto F_t(\bx)$ satisfies the proximal-PL inequality~\cite{karimi2016linear}, and an inexact gradient is available.  

We note that strong convexity implies the PL inequality. However, functions that satisfy the PL inequality are not necessarily convex; instead, they satisfy the notion of invexity~\cite{karimi2016linear}. Examples of cost functions that satisfy the PL inequality includes least squares (LS) and logistic regression, with applications that span learning and  feedback-based optimization. On the other hand, prime examples of costs that satisfy the proximal-PL condition are the LS with a sparsity-promoting regularize and an indicator function for a  polyhedral set (see, e.g.,~\cite{karimi2016linear} for additional examples).


The analysis is performed in terms of the instantaneous regret $r_t := F_t(\bx_t) - F_t^*$, where $F_t(\bx_t)$ is the cost achieved at time $t$ by the point $\bx_t$ produced by the algorithm and $F_t^*$ is the optimal value function (that one would have achieved if the problem~\eqref{eq:main-problem} was solved to convergence at time $t$). 

Motivating examples for considering stochastic gradient information are drawn from a variety of applications in learning and data-driven control; for example: i) settings where bandit and zeroth-order methods are utilized to estimate the gradient from (one or a few) functional evaluations~\cite{hajinezhad2017zeroth,tang2020distributed}; ii)~feedback-based optimization of networked systems, where errors in the gradient are due to measurement errors and asynchronous measurements~\cite{Bolognani_feedback_15,hauswirth2021optimization,ospina2022feedback}; and, iii)~online stochastic optimization settings~\cite{shames2020online,cao2020online}. 


\emph{Prior works}. Online (projected) gradient descent 
methods with exact gradient information have been investigated, and we refer to the 
representative works~\cite{selvaratnam2018Numerical,mokhtari2016online,madden2021bounds} as well as 
to references therein. A regret analysis was performed in, e.g.,~\cite{yi2016tracking,mokhtari2016online,chang2021online} (see also pertinent references therein), and the excess-risk was analyzed in~\cite{towfic2013distributed}. Inexact gradient information was considered in, e.g.,~\cite{besbes2015non,yi2016tracking}, where bounds in expectation on the regret incurred by the inexact online gradient descent were derived, and in~\cite{bedi2018tracking} where the distance from the unique trajectory of optimizers was bounded in expectation. Convergence results in expectation were provided in the context of online stochastic optimization in, e.g.,~\cite{shames2020online,cao2020online}. Convergence guarantees for  online 
stochastic gradient methods where drift and noise terms satisfy sub-Gaussian assumptions were provided
in~\cite{cutler2021stochastic}. Online projected gradient methods with sub-Weibull gradient error and a strongly convex cost are analyzed in~\cite{ospina2022feedback}.

We also acknowledge representative prior works on inexact and stochastic gradient and proximal-gradient methods for batch optimization in, e.g.,~\cite{schmidt2011convergence,devolder2014first,gannot2021frequency,rosasco2014convergence,atchade2017perturbed,moulines2011non,bertsekas1997gradient,li2020high} (see also references therein). In particular, almost sure convergence to a first-order stationary point is proved assuming only strong smoothness and a weak assumption on the noise in \cite{bertsekas1997gradient}; mean convergence under the PL inequality is shown in, e.g, \cite{khaled2020better}. High-probability convergence results assuming strong smoothness and norm sub-Gaussian noise were provided in e.g.,~\cite{li2020high}, and in~\cite{Harvey2} for strongly convex functions in the non-smooth setting. Finally, we also acknowledge prior works that investigate geometric conditions implying linear convergence of proximal gradient algorithms~\cite{bolte2007lojasiewicz,attouch2013convergence,bolte2010characterizations}. These works are for  static optimization.

\emph{Contributions}. We consider the cases \emph{c1)} and \emph{c2)} described above, and offer  the  following main contributions.   

\emph{(i)} We provide new bounds for the instantaneous regret $r_t$ in \emph{expectation} and in \emph{high probability} for the inexact online gradient descent, when the cost satisfies the PL inequality. The high-probability convergence results are derived by adopting a sub-Weibull model~\cite{vladimirova2020sub} for the gradient error. We also provide an \emph{almost sure}  result for the asymptotic behavior of the regret $r_t$. 

\emph{(ii)}  Similarly, we provide new bounds for the instantaneous regret $r_t$ in \emph{expectation} and in \emph{high probability} for the inexact online proximal-gradient descent method.

\emph{(iii)} For the case of static costs, our bounds provide contributions over~\cite{karimi2016linear,schmidt2011convergence,devolder2014first,gannot2021frequency,rosasco2014convergence,atchade2017perturbed,moulines2011non,bertsekas1997gradient,li2020high,bertsekas1997gradient,khaled2020better,li2020high} by considering a sub-Weibull model for the gradient error. In terms of bounds in expectation, this paper extends the results of in the context of static optimization to an online setting where the cost changes over time.  

To better highlight the merits of the bounds, is important to mention that the sub-Weibull distribution allows one to consider a variety of gradient error models in a unified manner; in fact, the sub-Weibull distribution includes sub-Gaussian distributions and sub-exponential distributions as sub-cases, as well as random variables whose probability density function has a finite support~\cite{vershynin_high-dimensional_2018}. The bounds we derived can be customized to sub-Gaussian  and sub-exponential distributions and for random variables with finite support by simply tuning the parameters of the sub-Weibull model. Furthermore,~\cite{bastianello2021stochastic} showed that intermittent updates can also be modeled using sub-Weibull random variables.

The rest of the paper is organized as follows. Section~\ref{sec:preliminaries} introduces relevant definitions and assumptions, and Section~\ref{sec:gradient} presents the main results for online gradient descent. Section~\ref{sec:prox-gradient} focuses on the online proximal-gradient method, and Section~\ref{sec:results} provides numerical results. Section~\ref{sec:conclusions} concludes the paper.

\section{Preliminaries}
\label{sec:preliminaries}

We start by introducing relevant definitions and assumptions that will be utilized throughout the paper\footnote{\emph{Notation}. Upper-case (lower-case) boldface letters will be used for matrices (column vectors); $(\cdot)^\top$ denotes transposition. For given column vectors $\bx, \by \in \mathbb{R}^n$, $\langle \bx,\by \rangle$ denotes the inner product and  $\|\bx\| := \sqrt{\langle \bx,\bx \rangle}$. Given a differentiable function $f: \cD \rightarrow \mathbb{R}$, defined over a domain $\cD \subseteq \mathbb{R}^n$ that is nonempty, $\nabla f (\bx)$ denotes the gradient of $f$ at $\bx$ (taken to be a column vector). 
$\mathcal{O}(\cdot)$ refers to the big-O notation, whereas $o(\cdot)$ refers to the little-o notation. For a given random variable $\xi \in \mathbb{R}$, $\mathbb{E}[\xi]$ denotes the expected value of $\xi$, and $\pr{\xi \leq \epsilon}$ denotes the probability of $\xi$ taking values smaller than or equal to $\epsilon$; furthermore, $\norm{\xi}_p := \mathbb{E}[|\xi|^p]^{1 / p}$, for any $p \geq 1$. Finally, $e$ will denote Euler's number.}. 

\subsection{Modeling and Definitions}

We consider functions $\{f_t\}_{t \in \mathbb{N}}$ and $\{g_t\}_{t \in \mathbb{N}}$, defined over an open ball $\mathcal{D} := \{ \bx \in \mathbb{R}^n: \|\bx\| < r$\} for some $r >0$, that satisfy the following assumptions. 

\vspace{.1cm}

\begin{assumption}
\label{as:f}
The function $\bx \mapsto f_t(\bx)$ is continuously differentiable and has a Lipschitz-continuous gradient over $\mathcal{D}$  for all $t$; i.e., $\exists~ L > 0$ such that $\| \nabla f_t(\bx) - \nabla f_t(\by) \| \leq L \|\bx - \by\|$ for any $\bx, \by \in \cD$, for all $t$. \QEDB
\end{assumption}

\vspace{.1cm}

\begin{assumption}
\label{as:g}
For every $t \in \mathbb{N}$, the function $\bx \mapsto g_t(\bx)$  is convex, proper, and lower semi-continuous, possibly non-differentiable over $\mathcal{D}$.   \QEDB
\end{assumption}

\vspace{.1cm}

Recall that the following inequality follows from the  Lipschitz-continuity of the gradient of $f_t$:
\begin{align}
\label{eq:smooth_def} 
f_t(\by) \leq f_t(\bx) + \langle \nabla f_t(\bx), \by - \bx \rangle + \frac{L}{2} \|\by - \bx\|^2 \, ,     
\end{align}
$\forall \, \bx, \by \in \cD$; this inequality will be utilized throughout the paper. Let $\cX^*_t := \arg \min_{\bx \in \mathbb{R}^n} F_t(\bx)
$, be the set of global minimizers of the problem~\eqref{eq:main-problem} at time $t$, and let $F_t^* := F_t(\bx_t^*)$, with $\bx_t^* \in \cX^*_t$. 
The following is assumed.

\vspace{.1cm}

\begin{assumption}
\label{as:bounded_optimal}
The set $\cX^*_t$ is non-empty for all $t$ and $\cX^*_t \subset \mathcal{D}$; furthermore,  $- \infty < F_t^*$  for all $t$. \QEDB
\end{assumption}

\vspace{.1cm}

The temporal variability of the problem~\eqref{eq:main-problem} could  be  measured  based on  how  fast  its  optimal solutions or optimal value functions change; see, for example,~\cite{besbes2015non,yi2016tracking,chang2021online} and references therein. More precisely, one can consider the change in the optimal value function as:
\begin{align}
\label{eq:variability_optimal}
\sigma_t := |F_{t}^* - F_{t-1}^*| \, . 
\end{align}
It will also be convenient to utilize the additional metrics $\tilde \phi_t(\bx)  := |F_{t}(\bx) - F_{t-1}(\bx)|$ and
\begin{align}
\label{eq:variability} 
\phi_t & := \sup_{\bx \in \mathcal{D}} \tilde \phi_t(\bx) \, .
\end{align}
For future developments, it is also convenient to define $\psi_t := \phi_t + \sigma_t$,  $\tilde \psi_t := \tilde \phi_t + \sigma_t$, and  $\bar{\psi}  := \sup_{t \in \mathbb{N}} \psi_t$. These metrics will be utilized to characterize the convergence of the online gradient and proximal-gradient methods. Whenever $g_t \equiv 0$ (this will be the main setting of Section~\ref{sec:gradient}), we will use the notation  $f_t^* = \min_{\bx \in \mathcal{D}} f_t(\bx)$ whenever convenient (in this case, it is clear that $F_t^* = f_t^*$); the definitions of $\sigma_t$, $\tilde \phi_t$, and $\phi_t$ remain unchanged. We also note that the case where $\sigma_t = 0$ and $\phi_t = 0$ for all $t$ corresponds to a static optimization problem (where the cost function does not change over time). 

We next recall the definition of the PL inequality and its generalization to composite cost functions~\cite{karimi2016linear}. 

\vspace{.1cm} 

\begin{definition}[Polyak-\L ojasiewicz (PL) Inequality]
\label{as:def_PL}
A continuously differentiable function $f:\mathcal{D}  \rightarrow \mathbb{R}$ satisfies the PL inequality over $\mathcal{D}$ if the following holds for some $\mu > 0$:
\begin{align}
\label{eq:pl} 
2 \mu (f(\bx) - f^*) \leq  \|\nabla f(\bx)\|^2, \,  \forall~\bx \in \mathcal{D} 
\end{align} 
where $f^*$ is the optimal value function. 
\end{definition}
\vspace{.1cm}

It is important to note that the PL inequality implies the quadratic bound $f(\bx) - f^* \geq \frac{\mu}{2} \|\bx - \bx^*\|^2$ for any global minimizer $\bx^*$.  As shown in~\cite{karimi2016linear}, strong convexity implies the PL inequality. However, functions that satisfy the PL inequality are not necessarily convex, instead, they 
satisfy the notion of invexity.  

\vspace{.1cm}

\begin{definition}[Proximal-PL Condition]
\label{as:def_proxPL}
Let $f:\mathcal{D}  \rightarrow \mathbb{R}$ be a continuously differentiable function and $g:\mathcal{D}  \rightarrow \mathbb{R}$ be a convex function. The function $F(\bx) := f(\bx) + g(\bx)$ satisfies the proximal-PL condition if the following holds:
\begin{align}
        \label{eq:prox_pl}
           2\mu (F(\bx)-F^*) \leq \mathcal{A}_{g} (\bx,\xi),
\end{align}
for all $\bx \in \cD$ and for some $\mu > 0$, where $\xi > 0$ and 
\begin{align}
\label{eq:defA}
     \mathcal{A}_{g} (\bx,\xi) := -\frac{2}{\xi} \min\limits_\by \Big\{ \langle \nabla f(\bx), \by-\bx \rangle + \frac{1}{2 \xi} \|\by-\bx\|^2 \notag \\
                +g(\by)-g(\bx) \Big\}.
 \end{align}
\end{definition}

\subsection{Sub-Weibull  random  variables}

In this section, we introduce the definition of sub-Weibull random variable (rv), which will be utilized to model the errors incurred by the inexact online gradient methods. 

\vspace{.1cm}

\begin{definition}[Sub-Weibull rv~\cite{vladimirova2020sub}]
\label{def:sub-weibull}
A random variable $X \in \mathbb{R}$ is sub-Weibull if  $\exists \, \theta > 0$ such that (s.t.) one of the following conditions is satisfied: 
\begin{enumerate}
	\item[(i)] $\exists \,\, K_1 > 0$ s.t. $\bbP[|X| \geq \epsilon] \leq 2 e^{- \left( \epsilon / K_1 \right)^{1 / \theta} }$, $\forall \, \epsilon~>~0$.
	\item[(ii)] $\exists \,\, K_2 > 0$ s.t. $\|X\|_k \leq K_2 k^\theta$, $\forall \, k \geq 1$. \hfill $\Box$
\end{enumerate}
\end{definition}
The parameters $K_1, K_2$ differ by a constant that depends on $\theta$. In particular, if (ii) holds with parameter $K_2$, then  (i) holds with $K_1 = \left( 2 e / \theta \right)^\theta K_2$. In this paper, we use the short-hand notation $X \sim \mathrm{subW}(\theta, K)$ to indicate that $X$ is a sub-Weibull rv according to Definition~\ref{def:sub-weibull}(ii). 

The coefficient $\theta$ is related to the rate of decay of the tails; in particular, the tails become heavier as the parameter $\theta$ grows larger. We note that the sub-Weibull class includes sub-Gaussian and sub-exponential rvs as sub-cases; in particular, if $\theta = 1/2$ and $\theta = 1$ we have sub-Gaussian and sub-exponential rvs, respectively. Furthermore, if a rv has a distribution with finite support, it belongs to the sub-Gaussian class (by Hoeffding’s inequality \cite[Theorem 2.2.6]{vershynin_high-dimensional_2018}) and, thus, to the sub-Weibull class.

\vspace{.1cm}

The following lemmas will be utilized throughout the paper to derive the main results.

\vspace{.1cm}

\begin{lemma} (\textit{Closure of sub-Weibull class}~\cite{bastianello2021stochastic}) Let $X_i \sim \mathrm{subW}(\theta_i, K_i)$, $i = 1,2$, based on Definition~\ref{def:sub-weibull}(ii). 
\begin{enumerate}	
	\item[(a)] \emph{Product by scalar:} Let $a \in \bbR$, then $a X_i \sim \mathrm{subW}(\theta_i, |a| K_i)$.
	\item[(b)] \emph{Sum by scalar:} Let $a \in \bbR$, then $a + X_i \sim \mathrm{subW}(\theta_i, |a| + K_i)$.
	\item[(c)] \emph{Sum:} Let $\{X_i, i = 1,2\}$ be possibly dependent; then, $X_1 + X_2 \sim \mathrm{subW}(\max\{ \theta_1, \theta_2 \}, K_1 + K_2)$.
\end{enumerate} \label{sec:closure}
\end{lemma}

\vspace{.1cm}

\begin{lemma} (\textit{Inclusion}~\cite{vladimirova2020sub}) 
Let $X \sim \mathrm{subW}(\theta, K)$ for some $\theta, K > 0$. Let $\theta', K'$ be s.t. $\theta' \geq \theta$, $K' \geq K$. Then, $X \sim \mathrm{subW}(\theta', K')$.  \hfill $\Box$ \label{sec:inclusion}
\end{lemma}

\vspace{.1cm}

\begin{lemma} (\textit{Powers of sub-Weibull rvs}~\cite{bastianello2021stochastic}) 
Let $X \sim \mathrm{subW}(\theta, K)$ for some $\theta, K > 0$, and let $a > 0$. Then, $X^a \sim \sw(a \theta, K^a \max\{ 1, a^{a \theta} \})$.  \hfill $\Box$ \label{sec:power}
\end{lemma}

\vspace{.1cm}


We  note that the definition of sub-Weibull rvs and their
properties do not require their mean to be zero. We conclude this section with the following high probability bound for a sub-Weibull rv. 

\vspace{.1cm}

\begin{lemma}[High probability bound]\label{lem:high-probability-bound}
Let $X \sim \mathrm{subW}(\theta, K)$ according to Definition \ref{def:sub-weibull}(ii), for some $\theta, K > 0$. Then, for any $\delta \in (0, 1)$, the bound:
\begin{equation}
| X | \leq \ K \log^\theta \left(2 \delta^{-1} \right) \left( \frac{2e}{\theta} \right)^\theta
\end{equation}
holds with probability $1 - \delta$. 
\hfill $\Box$
\end{lemma}

\section{Online Stochastic Gradient Descent}
\label{sec:gradient}

We start by considering the case where $g_t \equiv 0$ for all $t$; accordingly, the problem~\eqref{eq:main-problem} reduces here to: 
\begin{align}
\label{eq:problem-f}
\min_{\bx} f_t(\bx) 
\end{align}  
where we recall that $t \in \mathbb{N}$ is the time index. 
We consider the following \emph{inexact online gradient descent} (OGD): 
\begin{equation}
\label{eq:iogd} 
\bx_{t+1} = \bx_{t} - \eta \, \bv_t 
\end{equation}
where $\eta > 0$ is a given step-size, $\bv_t := \nabla f_t(\bx_{t}) + \be_t$ is the approximate gradient, and $\be_t$ is a stochastic error.  We are interested in studying the performance of~\eqref{eq:iogd} when the function $f_t$ satisfies the PL inequality~\eqref{eq:pl}, and the error $\|\be_t\|$ follows a sub-Weibull distribution. A  discussion on the sub-Weibull model as well as the PL inequality in the context of problems in learning and feedback-based optimization is provided in Section~\ref{sec:errormodel}. The main convergence results are presented next.

\subsection{Convergence in expectation and in high probability}
\label{sec:gradient_conv}

Since $g_t \equiv 0$, the instantaneous regret at time $t$ boils down here to $r_t = f_t(\bx_t) - f_t^*$. Throughout  this section, we assume that the gradient error has a sub-Weibull distribution, as formalized next. 

\vspace{.1cm}

\begin{assumption}[Sub-Weibull norm gradient error]
\label{as:sub-weibull}
The error is distributed as $\|\mathbf{e}_t\| \sim \sw(\theta, K_t)$, for some $\theta > 0$ and $K_t > 0$.
\end{assumption}

\vspace{.1cm}

We note that if each individual entry of the random vector $\mathbf{e}_t$ follows a sub-Weibull distribution, then $\|\be_t\|$  is a sub-Weibull rv. This can be proved by using~\cite[Lemma 3.4]{bastianello2021stochastic} and part (c) of Proposition~\ref{sec:closure}. In the following, we state the main results concerning the convergence of~\eqref{eq:iogd}.

\vspace{.1cm}

\begin{theorem}[Convergence of the stochastic OGD]
\label{thm:regret}
Let Assumptions~\ref{as:f},~\ref{as:bounded_optimal}, and~\ref{as:sub-weibull} hold, and assume that the map $\bx \mapsto f_t(\bx)$ satisfies the PL inequality, for some $\mu > 0$, for all $t$. Let $\{\bx_i\}_{i = 0}^t$  be a sequence generated by~\eqref{eq:iogd} with $\eta = 1/L$. The following bounds hold for \eqref{eq:iogd}:
\begin{enumerate}
    \item For all $t\in \mathbb{N}$:
    \begin{align}
    \label{eq:expected_regret}
     \hspace{-.4cm} \mathbb{E}[r_t]  \leq  \zeta^t r_0  + \sum_{\tau = 1}^t \zeta^{t - \tau} \left(\frac{1}{2L}  \mathbb{E}[\|\be_{\tau-1}\|^2]  + \psi_\tau \right) 
     \end{align}
    where $\zeta := (1-\frac{\mu}{L})$. 
    
    \item For any $\delta\in(0,1)$, then the following bound holds with probability $1-\delta$:
    \begin{align}
    \label{eq:high_regret}
     \hspace{-.4cm} r_t \leq h(\theta, \delta) \left(\zeta^t r_0 + \sum_{\tau = 1}^t \zeta^{t - \tau} \left( \frac{4^\theta}{2L} K_{\tau-1}^2 + \psi_\tau  \right) \right)
     \end{align}
     where 
     $ h(\theta, \delta) := \log^{2 \theta}(2 \delta^{-1}) \left(\frac{e}{\theta} \right)^{2 \theta}$.
\end{enumerate}

\end{theorem}

\vspace{.1cm}

\begin{corollary}[Asymptotic convergence]
\label{cor:asymp_regret}
Under the same assumptions of Theorem~\ref{thm:regret}, it holds that
    \begin{align}
     \limsup_{t \rightarrow \infty} r_t \leq  \frac{1}{2 \mu} \bar{e} + \frac{L}{\mu} \bar \psi \,\,\,\, \textrm{a.s.}
     \end{align}
where $\bar{e} = \sup_t \{\mathbb{E}[\|\be_{t}\|^2]\}$. 
\end{corollary}

\vspace{.1cm}

Before providing examples of applications and the proof of the results, some remarks are in order. 

\vspace{.1cm}

\begin{remark}[Static optimization~\cite{karimi2016linear}]
\label{rem:static}
When the optimization problem~\eqref{eq:problem-f} is time-invariant (i.e., $f_t(x) = f(x)$ for all $t \in \mathbb{N}$), then,~\eqref{eq:expected_regret} is similar to~\cite[Thm.~4]{karimi2016linear} (where a different step-size was used). However, relative to~\cite{karimi2016linear}, we provide the following bound in high probability
\begin{align}
    \label{eq:high_regret2}
     \hspace{-.4cm} r_t \leq   h(\theta, \delta) \left(\zeta^t r_0 + \sum_{\tau = 1}^t \zeta^{t - \tau} \frac{4^\theta}{2L} K_{\tau-1}^2  \right)
\end{align}
which holds with probability $1 - \delta$ for any $\delta \in (0,1)$; this bound can be derived from~\eqref{eq:high_regret} by setting $\psi_\tau = 0$ for all $\tau = 1, \ldots, t$. 
 \hfill \QEDB 
\end{remark}

\vspace{.1cm}

\begin{remark}[Alternative bound in expectation]
\label{rem:tighterBound}
An alternative bound in expectation can be expressed as 
\begin{align}
    \label{eq:expected_regret2}
     \mathbb{E}[r_t]  \leq  \zeta^t r_0  + \sum_{\tau = 1}^t \zeta^{t - \tau} \left(\frac{1}{2L}  \mathbb{E}[\|\be_{\tau-1}\|^2]  + \mathbb{E}[\tilde{\psi}_\tau] \right) 
     \end{align}
where $\tilde{\psi}_\tau = \sigma_\tau + \tilde{\phi}_\tau$, and $\tilde{\phi}_\tau :=  |F_{\tau}(\bx_\tau) - F_{\tau-1}(\bx_\tau)|$ (where the expectation $\mathbb{E}[\tilde{\psi}_\tau] $ is taken with respect to the error $\be_{\tau-1}$, conditioned on a filtration). This leads to a tighter bound relative to~\eqref{eq:expected_regret}.      
\hfill \QEDB   
\end{remark}

\vspace{.1cm}

\begin{remark}[Markov's inequality] 
An alternative high probability bound 
can be obtained by using~\eqref{eq:expected_regret} and Markov's inequality. However, the resulting bound would have a dependence $\delta^{-1}$; on the other hand, our bound has a $\log(\delta^{-1})$ dependence on $\delta$. 
 \hfill \QEDB 
\end{remark}

\subsection{Remarks on applications and error model}
\label{sec:errormodel}

In this section, we provide some examples of applications that are relevant to our setting. 

\vspace{.1cm}
\begin{example}[Online least-squares]
A function $f_t(\bx) = h_t(\bA_t \bx)$, with $h_t: \mathbb{R}^d \rightarrow \mathbb{R}$ a $\nu$-strongly convex function and $\bA_t \in \mathbb{R}^{d \times n}$ satisfies the PL inequality~\cite{karimi2016linear}. This class includes the least-squares (LS) problem by setting $f_t(\bx) = \frac{1}{2} \|\bA_t \bx - \bb_t\|^2$. Note that, when the matrix $\bA_t$ is not full-column rank, one can utilize the results of this paper to establish linear convergence of OGD for the under-determined LS problem.  
\end{example}

\vspace{.1cm}

\begin{example}[Online Logistic regression]
 The logistic regression cost $f_t(\bx) = \sum_{i = 1}^d \log(1 + \exp(b_{i,t} \ba_{i,t}^\top \bx))$, with $b_{i,t} \in \mathbb{R}$ and $\ba_{i,t} \in \mathbb{R}^n$,  satisfies the PL inequality~\cite{karimi2016linear}. 
\end{example}

\vspace{.1cm}

\begin{example}[Optimization of LTI systems] 
\label{ex:feedbackoptimization}
Consider the algebraic representation of a stable linear time-invariant system $\by_t = \bG \bx + \bH \bw_t$, where $\bx$ is the vector of controllable inputs and $\bw_t$ are unknown exogenous disturbances. Suppose that $f_t(\bx) = \frac{1}{2}\|\bG \bx + \bH \bw_t - \bar{\by}_t\|^2$ with $\bar{\by}_t$ a time-varying reference. Since $\bw_t$ is unknown, one way to compute $\nabla f_t(\bu_t)$ is $\bv_t = \bG^\top (\hat{\by}_t - \bar{\by}_t)$, where $\hat{\by}_t$ is a (noisy) measurement of the output $\by_t$~\cite{Bolognani_feedback_15,ospina2022feedback}. 
\end{example}

\vspace{.1cm}

\begin{example}[Training of neural networks]
We  refer the reader to recent discussions on the PL inequality in the context of training of neural networks in, e.g.,~\cite{li2017convergence}. The proposed framework may capture the case where stochastic gradient methods are utilized to train a neural network in an online fashion. 
\end{example}

\vspace{.1cm}

In terms of gradient information, the error $\be_t$ may arise in the following (application-specific) scenarios: 

\noindent \emph{(i)} A subset of the data points available at time $t$ are utilized to compute the gradient; for instance, in the Examples 1-2, one may utilize the data points $\{\ba_{i,t}, b_{i,t}\}_{i \in \cS_t}$, with $|\cS_t| < d$.   

\noindent \emph{(ii)} Bandit and zeroth-order methods are utilized to estimate the gradient~\cite{hajinezhad2017zeroth,tang2020distributed}. 

\noindent \emph{(iii)} In an online stochastic optimization setting~\cite{shames2020online}, i.e. when $f_t(\bx) = \ev{\ell_t(\bx,\bz)}$ for a given loss $\ell_t: \mathbb{R}^n \times \mathbb{R}^d \rightarrow \mathbb{R}$ and a random variable $\bz$, the approximate gradient  $\bv_t$ may be computed using a single sample or a mini-batch. 

\noindent \emph{(iv)} In measurement-based algorithms as in Example~\ref{ex:feedbackoptimization}, measurement errors and asynchronous measurements render the computation of the gradient inexact.

\subsection{Proofs}
\label{sec:proof}

In this section, we present the proof of  Theorem~\ref{thm:regret}.

We start by using~\eqref{eq:smooth_def} with $\by = \bx_{t+1}$ and $\bx = \bx_{t}$, where $\bx_{t}$ and $\bx_{t+1}$ are generated by~\eqref{eq:iogd}; this allows us to obtain:
\begin{subequations}
\begin{align}
    f_{t}(\bx_{t+1}) & \leq f_t(\bx_t) -\frac{1}{L}  \langle \nabla f_t(\bx_t),\bv_t  \rangle + \frac{1}{2L} \| \bv_t \|^2 \\
    &=f_t(\bx_t) -\frac{1}{L}  \langle \nabla f_t(\bx_t), \nabla f_t(\bx_t) +  \mathbf{e}_t \rangle \notag \\
    &\quad+ \frac{1}{2L} \|\nabla f_t(\bx_t)+\mathbf{e}_t \|^2 \\
    &=f_t(\bx_t)-\frac{1}{L} \|\nabla f_t(\bx_t)\|^2 -\frac{1}{L} \langle \nabla f_t(\bx_t), \mathbf{e}_t \rangle \notag \\
    &\quad+ \frac{1}{2L} \|\nabla f_t(\bx_t)+\mathbf{e}_t\|^2 \\
    &=f_t(\bx_t)-\frac{1}{L} \|\nabla f_t(\bx_t)\|^2 -\frac{1}{L} \langle \nabla f_t(\bx_t), \mathbf{e}_t \rangle \notag \\
    &\quad+ \frac{1}{2L} \langle \nabla f_t(\bx_t)+\mathbf{e}_t, \nabla f_t(\bx_t)+\mathbf{e}_t \rangle \\
    &=f_t(\bx_t)-\frac{1}{L} \|\nabla f_t(\bx_t)\|^2 -\frac{1}{L} \langle \nabla f_t(\bx_t), \mathbf{e}_t \rangle \notag \\
    &\quad+ \frac{1}{2L}(\| \nabla f_t(\bx_t)\|^2+ 2 \langle \nabla f_t(\bx_t),\mathbf{e}_t \rangle +\|\mathbf{e}_t\|^2)\\
   &=f_t(\bx_t) - \frac{1}{2L} \|\nabla f_t(\bx_t)\|^2+ \frac{1}{2L} \|\mathbf{e}_t\|^2 \, .
\end{align}
\end{subequations}
Next, adding $- f_t^*$ on both sides and using the PL inequality~\eqref{eq:pl}, one gets:
\begin{align}
    &f_t(\bx_{t+1})-f_t^* \leq - \frac{1}{2L} \|\nabla f_t(\bx_t)\|^2 + \frac{1}{2L} \|\mathbf{e}_t\|^2+f_t(\bx_t)-f_t^* \notag \\
    & \leq -\frac{\mu}{L}(f_t(\bx_t)-f_t^*)+f_t(\bx_t)-f_t^* + \frac{1}{2L}\|\mathbf{e}_t\|^2 \, .
\end{align}
Next, adding $-f_{t+1}^*$ and $f_{t+1}(\bx_{t+1})$ on both sides, using the definition of regret $r_t$, and applying the definitions of $\zeta$ and $\psi_t$ we obtain the stochastic inequality 
$$
r_{t} \leq \zeta \, r_{t-1} +  \frac{1}{2L} \|\be_{t-1}\|^2 + \psi_t , $$ 
which holds almost surely. Unraveling, we get 
\begin{align}
\label{eq:randombound}    
    r_t \leq \kappa_t +\frac{1}{2L}\sum_{i=1}^{t}\zeta^{t-i}\|\mathbf{e}_{i-1}\|^2 
\end{align}
where $\kappa_t := \zeta^t r_0 + \sum_{i=1}^{t} \zeta^{t-i} \psi_i$ for brevity. Taking the expectation on both sides, we get~\eqref{eq:expected_regret}. 

For the high-probability bound~\eqref{eq:high_regret}, 
recall that $\|\mathbf{e}_i\| \sim \sw(\theta, K_i)$; by Lemma~\ref{sec:power}, setting $a = 2$ we get that $\|\mathbf{e}_i\|^2$ is a sub-Weibull rv and, in particular, $\|\mathbf{e}_i\|^2 \sim \sw(2 \theta, 4^\theta K_i^2)$. Next, using the closure properties (a), (b), and (c) in Lemma~\ref{sec:closure} and the fact that $\zeta > 0$, we have that the right-hand-side of~\eqref{eq:randombound}  is a sub-Weibull rv; in particular, 
\begin{align}
\label{eq:errortotal}  
    \kappa_t + \frac{1}{2L} \sum_{i=1}^{t}\zeta^{t-i}\|\mathbf{e}_{i-1}\|^2  \sim \sw \left(2 \theta, K^\prime \right) \, .
\end{align}
where 
$$
K^\prime = \kappa_t + \frac{4^\theta}{2L} \sum_{i=1}^{t}\zeta^{t-i} K_{i-1}^2 
$$
Using Lemma~\ref{lem:high-probability-bound}, the high-probability bound~\eqref{eq:high_regret} follows. 

The proof of Corollary~\ref{cor:asymp_regret} follows similar steps as in~\cite[Corollary~4.8]{bastianello2021stochastic}, and  is omitted.

\section{Stochastic Proximal-Gradient Method}
\label{sec:prox-gradient}

We now turn the attention to the time-varying problem~\eqref{eq:main-problem}, with the cost satisfying   the Assumptions~\ref{as:f}-\ref{as:bounded_optimal}. Throughout this section, we further assume that the cost function $F_t(\bx)$ satisfies the proximal-PL inequality~\eqref{eq:prox_pl}, for a given $\mu > 0$. As discussed in~\cite{karimi2016linear}, an important example of cost satisfying the proximal-PL inequality is the $\ell_1$-regularized least squares problem; additional examples of costs include (see the discussion in~\cite[Appendix F]{karimi2016linear}): 

\begin{enumerate}
    \item $F_t(\bx) = f_t(\bA \bx) + g_t(\bx)$, with $f_t$ strongly convex, $g_t$ the indicator function for a polyhedral set, and $\bA$ a given matrix. 

\item The case where $f_t$ is convex, and $F_t$ satisfies the quadratic growth condition.  

\item The case where $F_t$ satisfies the Kurdyka-\L ojasiewicz inequality or the proximal exponential bound. 
\end{enumerate} 

Consider then the stochastic online proximal-gradient method (OPGM), which involves the following step:
\begin{align}
    \label{eq:opgm}
        \bx_{t+1} = \mathrm{prox}_{\frac{1}{L} g_t} \left\{\bx_t - \frac{1}{L} \bv_t \right\} \, , \,\,\,\, t \in \mathbb{N}
\end{align}
where $\bv_t$ is again an estimate of $\nabla f_t(\bx_t)$,  $\mathrm{prox}_{\frac{1}{L} g_t}: \mathbb{R}^n \rightarrow \mathbb{R}^n$ denotes the proximal operator, and the step-size is taken to be $1/L$. 

We are now interested in analyzing the behavior of~\eqref{eq:opgm} in terms of regret $r_t = F_t(\bx_t) - F_t^*$, where we recall that $F_t^*$ is the optimal value function, when the function $F_t$ satisfies proximal-PL inequality and the error $\|\be_t\|$ follows a sub-Weibull distribution. The main convergence result for~\eqref{eq:opgm} is  stated next.

\vspace{.1cm}

\begin{theorem}[Convergence of the stochastic OPGM]
\label{thm:regretOPGM}
Let Assumptions~\ref{as:f}--\ref{as:sub-weibull} hold. Assume further that the function $\bx \mapsto F_t(\bx)$ satisfies the proximal-PL inequality for some $\mu > 0$, for all $t$. Let $\{\bx_i\}_{i = 0}^t$  be a sequence generated by~\eqref{eq:opgm}. Then:
\begin{enumerate}
    \item For all $t\in \mathbb{N}$:
    \begin{align}
    \label{eq:expected_regretOPGM}
     \hspace{-.5cm} \mathbb{E}[r_t]  \leq  \zeta^t r_0  + \sum_{\tau = 1}^t \zeta^{t - \tau} \left( 2 D \mathbb{E}[\|\be_{\tau-1}\|] + \psi_\tau \right) 
     \end{align}
    where $\zeta = (1-\frac{\mu}{L})$ and $D$ is the diameter of $\mathcal{D}$. 
    
    \item If $\delta\in(0,1)$, then with probability $1-\delta$: 
    \begin{align}
    \label{eq:high_regretOPGM}
     \hspace{-.5cm} r_t \leq h_p(\theta, \delta) \left(\zeta^t r_0 + \sum_{\tau = 1}^t \zeta^{t - \tau} \left( 2D K_{\tau-1} + \psi_\tau  \right) \right)
     \end{align}
     where 
     $ h_p(\theta, \delta) := \log^{\theta}(2 \delta^{-1}) \left(\frac{2e}{\theta} \right)^{\theta}$.
\end{enumerate}

\end{theorem}

\vspace{.1cm}

A result for the asymptotic convergence of the OPGM similar to Corollary~\ref{cor:asymp_regret} can be derived too, but it is omitted to avoid repetitive arguments. Similar considerations as in Remark~\ref{rem:static} can also be drawn. 

We note that when $g_t$ is the indicator function for a bounded polyhedron, the constant $D$ can be replaced by the diameter of the polyhedron. We also note that tighter bounds could be derived by introducing a filtered probability space (this will be clear in the proof, where the inner product between the iterates and the gradient error appears); however, this is left as a future extension.

To outline the proof of the theorem, we first note that step~\eqref{eq:opgm} is equivalent to 
    \begin{align}
    \label{eq:opgm2}
        \bx_{t+1} =\arg \min_{\bx } \langle \bv_t, \bx -\bx_t \rangle & + \frac{L}{2} \|\bx -\bx_t\|^2  +g_t(\bx)-g_t(\bx_t) .
    \end{align}
We also recall the definition of $\mathcal{A}_{g_t}(\bx_t,1/L)$ in~\eqref{eq:defA}, and define  $\tilde{\mathcal{A}}_{g_t}(\bx_t,1/L)$ as:
    \begin{align}
        \label{eq:tildeA}
        \tilde{\mathcal{A}}_{g_t}(\bx_t,1/L)&:=-2L\min\limits_{\by} \{\langle \bv_t,\by-\bx_t \rangle+\frac{L}{2} \|\by-\bx_t\|^2 \nonumber \\
        & \hspace{2.8cm} +g_t(\by)-g_t(\bx_t) \} \, .
    \end{align}
Lastly, for any $\bx \in \mathcal{D}$, we recall that  $\mathbf{e}_t \in \mathbb{R}^n$ is the gradient error, i.e.,  $\bv_t=\nabla f_t(\bx)+\mathbf{e}_t$, and $\|\be_t\| \sim \sw(\theta, K_t)$.  

\vspace{.1cm}

\emph{Proof of Theorem~\ref{thm:regretOPGM}.} We start by recalling that $F_{t+1}(\bx_{t+1})=f_{t+1}(\bx_{t+1})+g_{t+1}(\bx_{t+1})$; adding and subtracting $F_t(\bx_{t+1})$ on the right-hand-side, and using the definition~\eqref{eq:variability}, we get   
    \begin{align}
        & F_{t+1}(\bx_{t+1}) \leq \phi_{t+1} +  f_t(\bx_{t+1})  +g_t(\bx_{t+1})+g_t(\bx_t)-g_t(\bx_t) \nonumber \\
        &\leq  \phi_{t+1}+ f_t(\bx_t)+\langle \nabla f_t(\bx_t),\bx_{t+1}-\bx_t \rangle \nonumber \\
        &\quad+\frac{L}{2} \|\bx_{t+1}-\bx_t\|^2 +g_t(\bx_{t+1})+g_t(\bx_t)-g_t(\bx_t) 
   \end{align}     
where we used~\eqref{eq:smooth_def} in the last step. Next, we add and subtract $\be_t$ in the inner product and use the definition $\bv_t=\nabla f_t(\bx)+\mathbf{e}_t$ to obtain: 
 \begin{align}       
         F_{t+1}(\bx_{t+1}) & \nonumber \\
         & \hspace{-1.1cm} \leq F_t(\bx_t)+ \langle  \bv_t, \bx_{t+1}-\bx_t \rangle+\frac{L}{2} \|\bx_{t+1}-\bx_t\|^2 \nonumber \\
        & \hspace{-.7cm} +g_t(\bx_{t+1})-g_t(\bx_t) -\langle \mathbf{e}_t,\bx_{t+1}-\bx_t \rangle + \phi_{t+1} \\
        & \hspace{-1.1cm} \leq F_t(\bx_t)-\frac{1}{2L} \tilde{\mathcal{A}}_{g_t}(\bx_t,1/L)+\phi_{t+1} -\langle \mathbf{e}_t,\bx_{t+1}-\bx_t \rangle 
\end{align}
where we used the definition~\eqref{eq:tildeA}. Adding and subtracting $\mathcal{A}_{g_t}(\bx_t,1/L)$ on the right-hand-side, we get      
    \begin{align}    
        F_{t+1}(\bx_{t+1}) &\leq F_t(\bx_t)+\frac{1}{2L} |\mathcal{A}_{g_t}(\bx_t,1/L) - \tilde{\mathcal{A}}_{g_t}(\bx_t,1/L)| \nonumber \\
        & \hspace{-.8cm} -\frac{1}{2L} \mathcal{A}_{g_t}(\bx_t,1/L) +\phi_{t+1}-\langle \mathbf{e}_t,\bx_{t+1}-\bx_t \rangle \, .
    \end{align}
Let $\varepsilon_t := |\mathcal{A}_{g_t}(\bx_t,1/L) - \tilde{\mathcal{A}}_{g_t}(\bx_t,1/L)|$ for brevity; using the definition of $\mathcal{A}_{g_t}(\bx_t,1/L)$ in~\eqref{eq:defA}  and  subtracting $F_{t+1}^*$ on both sides, we get: 
\begin{subequations}
    \begin{align}
        &F_{t+1}(\bx_{t+1})-F_{t+1}^* \leq 
        F_t(\bx_t) - F_{t+1}^* -\frac{\mu}{L}(F_t(\bx_t)-F_t^*) \nonumber  \\
        &\quad +\frac{1}{2L} \varepsilon_t  -\langle \mathbf{e}_t,\bx_{t+1}-\bx_t \rangle +\phi_{t+1}  \\
        & \leq 
        F_t(\bx_t) - F_{t}^* -\frac{\mu}{L}(F_t(\bx_t)-F_t^*) +\frac{1}{2L} \varepsilon_t  -\langle \mathbf{e}_t,\bx_{t+1}-\bx_t \rangle \nonumber  \\
        &\quad  +\phi_{t+1} + \sigma_{t+1} \\
        & = (1-\frac{\mu}{L})(F_t(\bx_t)-F_t^*)+\psi_{t+1}+\frac{1}{2L} \varepsilon_t -\langle \mathbf{e}_t,\bx_{t+1}-\bx_t \rangle  \\
        & = (1-\frac{\mu}{L})(F_t(\bx_t)-F_t^*)+\psi_{t+1}+\frac{1}{2L} \varepsilon_t + \|\mathbf{e}_t\| D  \label{eq:boundF}
    \end{align}
\end{subequations}
where we used the definition of $\sigma_{t+1}$ and $D$ is the diameter of $\mathcal{D}$.  

We now bound $\frac{1}{2L} \varepsilon_t$; from the definitions~\eqref{eq:defA} and~\eqref{eq:tildeA}, we have that 
\begin{align}
\frac{1}{2L} \varepsilon_t  = &  \Big| \min_{\by} \left\{\langle \nabla f_t(\bx_t) , \by - \bx_t \rangle +\frac{L}{2} \|\by-\bx_t\|^2 +g_t(\by) \right\} \nonumber \\
 & \hspace{.2cm} -  \min_{\bz} \left\{\langle \bv_t , \bz - \bx_t \rangle +\frac{L}{2} \|\bz-\bx_t\|^2 +g_t(\bz) \right\} \Big| .
\end{align}
From~\eqref{eq:opgm2}, one can notice that the minimizer of $\langle \bv_t , \bz - \bx_t \rangle +\frac{L}{2} \|\bz-\bx_t\|^2 +g_t(\bz)$ is $\bx_{t+1}$ (the constant term $g_t(\bx_{t})$ does not modify the minimizer); thus, substituting $\bz$ with $\bx_{t+1}$ we get
\begin{align}
\frac{1}{2L} \varepsilon_t  = &  \Big| \min_{\by} \left\{\langle \nabla f_t(\bx_t) , \by - \bx_t \rangle +\frac{L}{2} \|\by-\bx_t\|^2 +g_t(\by) \right\} \nonumber \\
 & \hspace{-.7cm} - \langle \bv_t , \bx_{t+1} - \bx_t \rangle - \frac{L}{2} \|\bx_{t+1}-\bx_t\|^2 - g_t(\bx_{t+1})  \Big|. \hspace{-.1cm}
\end{align}
Next, one has that $\min_{\by} \{\langle \nabla f_t(\bx_t) , \by - \bx_t \rangle +\frac{L}{2} \|\by-\bx_t\|^2 +g_t(\by) \} \leq \langle \nabla f_t(\bx_t) , \by - \bx_t \rangle +\frac{L}{2} \|\by-\bx_t\|^2 +g_t(\by)$ for any $\by \in \mathcal{D}$, and, thus 
\begin{align}
\frac{1}{2L} \varepsilon_t  \leq &  \Big| \langle \nabla f_t(\bx_t) , \by - \bx_t \rangle +\frac{L}{2} \|\by-\bx_t\|^2 +g_t(\by)  \nonumber \\
 & \hspace{-.7cm} - \langle \bv_t , \bx_{t+1} - \bx_t \rangle - \frac{L}{2} \|\bx_{t+1}-\bx_t\|^2 - g_t(\bx_{t+1})  \Big| 
\end{align}
for any $\by \in \mathcal{D}$. Pick $\by = \bx_{t+1}$; then, we have that:
\begin{subequations}
\begin{align}
\frac{1}{2L} \varepsilon_t  \leq &  \Big| \langle \nabla f_t(\bx_t) , \bx_{t+1} - \bx_t \rangle - \langle \bv_t , \bx_{t+1} - \bx_t \rangle \Big| \\
= &  \Big|  \langle \be_t , \bx_t - \bx_{t+1}  \rangle  \Big| \leq  \|\mathbf{e}_t\| D .
\end{align}
\end{subequations}

Therefore, letting $r_t = F_t(\bx_{t})-F_{t}^*$ for brevity,  we get the stochastic recursion
\begin{align}
        r_{t+1} \leq &  (1-\frac{\mu}{L}) r_t+\psi_{t+1}  + 2 D   \|\mathbf{e}_t\|   \label{eq:boundF2}
    \end{align}
which holds almost surely. By applying recursively~\eqref{eq:boundF2} from $\tau = 0$ to $\tau = t$, we get
    \begin{align}
    \label{eq:boundF3}
     \hspace{-.2cm} r_t  \leq  \zeta^t r_0  + \sum_{i = 1}^t \zeta^{t - i} \left( 2 D  \|\mathbf{e}_{i-1}\| + \psi_i \right) \, .
     \end{align}
Taking the expectation on both sides of~\eqref{eq:boundF3}, the bound~\eqref{eq:expected_regretOPGM} follows.

To show~\eqref{eq:high_regretOPGM}, recall first that $\|\be_t\| \sim \sw(\theta, K_t)$, and let  $\kappa_t := \zeta^t r_0 + \sum_{i=1}^{t} \zeta^{t-i} \psi_i$ for brevity so that~\eqref{eq:boundF3} can be rewritten as $r_t  \leq \kappa_t   + 2D \sum_{i = 1}^t \zeta^{t - i}  \|\mathbf{e}_{i-1}\|$. Using the closure properties (a), (b), and (c) in Lemma~\ref{sec:closure} and the fact that $\zeta > 0$ and $\kappa_t \geq 0$, we have that the right-hand-side of this inequality  is a sub-Weibull rv; in fact, 
\begin{align}
\label{eq:errortotal}  
    \kappa_t   + 2D \sum_{i = 1}^t \zeta^{t - i}  \|\mathbf{e}_{i-1}\| \sim \sw \left(\theta, K^{''} \right) \, .
\end{align}
where $K^{''} = \kappa_t + 2D \sum_{i=1}^{t}\zeta^{t-i} K_{i-1}
$. 
Using Lemma~\ref{lem:high-probability-bound}, the high-probability bound~\eqref{eq:high_regret} follows.

\section{Illustrative Numerical Results}
\label{sec:results}

We provide two illustrative numerical experiments. The first one is based on a time-varying LS regression problem; then, we consider a problem related to real-time demand response in power grids. 

\textbf{Least-squares problem}. We consider a time-varying LS regression problem, with the following cost at time $t$:
\begin{align}
    f_t(\bx) =  \frac{1}{2}\|\bA \bx - \bb_t\|^2
\end{align}
where $\bA \in \mathbb{R}^{d \times n}$ and $\bb_t \in \mathbb{R}^d$; this cost satisfies the PL inequality, as shown in~\cite{karimi2016linear}.

We consider the case $n=10$ and $d=20$.  The matrix $\bA$ is generated by defining its singular value decomposition; for its left and right-singular vectors, we sampled two orthogonal matrices, $\bU\in\mathbb{R}^{d\times d}$ and $\bV\in\mathbb{R}^{n\times n}$, and we let its singular values be equally spaced from $\mu = 0.1$ to $L = 1$. We generated $\bb_t$ as $\bb_t = \bA \bx_t^* + \br_t$, where  the optimal parameter $\bx_t^*$ evolves via a random walk, i.e., $\bx_t^* = \bx_{t-1}^* + \bs$ with $\bs \sim \mathcal{N}(\mathbf{0}, 0.1 \bI)$, and $\br_t$ is a Gaussian vector  $\mathcal{N}(\mathbf{0}, 10^{-3} \bI)$ (we set $\bx_{0}^*$ to the vector of all ones). 

\begin{figure}[t]
    \centering
    \includegraphics[width=\columnwidth]{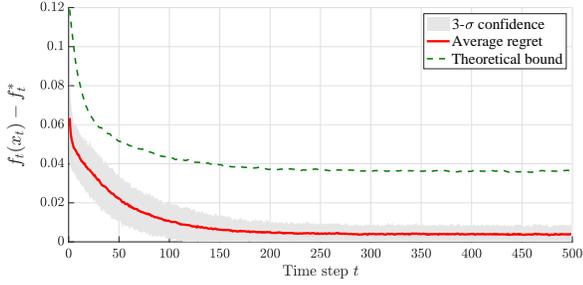}
    \caption{Inexact OGD: Evolution of average regret obtained experimentally, the empirical $3$-standard deviation confidence interval, and the theoretical bound.}
    \label{fig:regret}
\end{figure}

We  corrupt the gradient with a random vector $\be_t$, which is modelled as a Gaussian vector  $\mathcal{N}(\mathbf{0}, 10^{-3} \bI)$; we note that, if $\be_t$ is a Gaussian vector, then $\|\be_t\|^2$ is a sub-Weibull random variable. The regret is computed using a Monte Carlo approach, with 100 tests. Accordingly, Figure~\ref{fig:regret} illustrates the evolution of the expected regret obtained by averaging the trajectories of the instantaneous regret $r_t$ over the various runs, the empirical $3-\sigma$ confidence interval, and the theoretical bound~\eqref{eq:expected_regret2}. The figure validates the convergence results for the inexact OGD and, since $\bb_t$ continuously changes, the average $r_t$ exhibits a plateau.

\textbf{Real-time demand response problem}.
We consider an example in the context of a power distribution grid serving residential houses or commercial facilities. We consider $n$ controllable distributed energy resources (DERs) providing services to the main grid; precisely, consider the setting where the vector  $\bx$ collects the active power outputs of the DERs, and assume the algebraic relationship $p_{0,t} = \ba_x^\top \bx + \ba_w^\top \bw_t$ for the net active power  at the point of common coupling, where $\ba_x \in \mathbb{R}^n$ and $\ba_w \in \mathbb{R}^w$ are sensitivity coefficients, and $\bw_t \in \mathbb{R}^w$ is a vector collecting active powers of uncontrollable devices; in particular, $\ba_x$ and $\ba_w$ can be set to the vector of all ones when line losses are negligible, or they are derived based on a linearized model for the power flow equations in case of resistive lines~\cite{Bolognani_feedback_15}.  Consider the following time-varying optimization problem for real-time management of DERs:
\begin{equation}
\label{eq:tip_sg}
\min_{\bx} \, \frac{1}{2} \left(\ba_x^\top \bx + \ba_w^\top \bw_t - p_{0,t}^{\mathrm{ref}} \right)^2 + \mathbb{I}_{\{\bB \bx \leq \bc\}}
\end{equation}
where $p_{0,t}^{\mathrm{ref}}$ is a time-varying reference point for the net active power  at the point of common coupling $p_{0,t}$, and $\mathbb{I}_{\{\bB \bx \leq \bc\}}$ is the set indicator function for the set $\{\bx \in \mathbb{R}^n: \bB \bx \leq \bc\}$ modeling box constraints for the active powers. For example, $p_{0,t}$ may be an automatic control generation (ACG) signal, a flexible ramping signal,  or a demand response setpoint.  We note that the cost~\eqref{eq:tip_sg} satisfies the proximal-PL inequality~\cite{karimi2016linear}. The main challenge behind applying a proximal-gradient descent to~\eqref{eq:tip_sg} is that the vector $\bw_t $ is unknown; we therefore consider the approach of, e.g.,~\cite{Bolognani_feedback_15}, where measurements of $p_{0,t}$ are utilized to estimate the gradient in lieu of the model $\ba_x^\top \bx + \ba_w^\top \bw_t $. Precisely, we compute the approximate gradient as
\begin{equation}
\bv_t = \ba_x (\hat{p}_{0,t} - p_{0,t}^{\mathrm{ref}})
\end{equation}
where $\hat{p}_{0,t}$ is a measurement of $p_{0,t}$ collected at time $t$. Since measurements of $p_{0,t}$  may be affected by errors or by outliers, $\bv_t$ does not in general coincide with the true gradient $\ba_x (\ba_x^\top \bx_t + \ba_w^\top \bw_t  - p_{0,t}^{\mathrm{ref}})$.

\begin{figure}[t]
    \centering
    \includegraphics[width=\columnwidth]{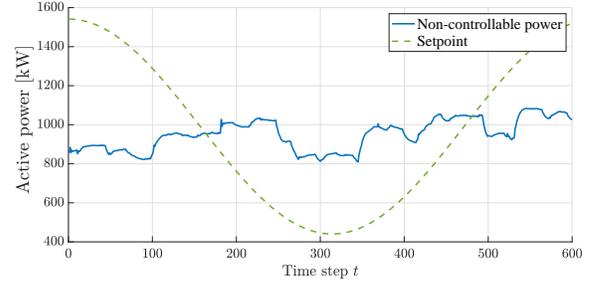}
    \caption{Demand response application: non-controllable power $\ba_w^\top \bw_t$ and reference point  $p_{0,t}^{\mathrm{ref}}$  for the  active power  $p_{0,t}$. }
    \label{fig:regret_data}
\end{figure}
\begin{figure}[t]
    \centering
    \includegraphics[width=\columnwidth]{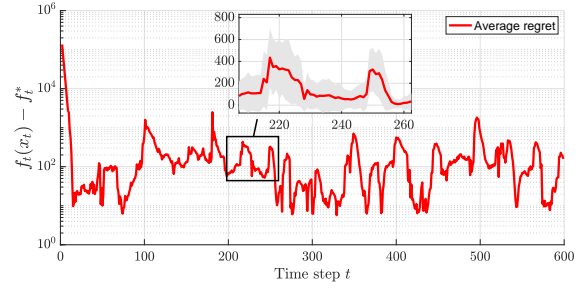}
    \caption{Demand response application: Evolution of average regret obtained experimentally; the zoomed area also provides the empirical $3-\sigma$ confidence interval.}
    \label{fig:regret_demand}
\end{figure}

As an example, we consider the case where $N = 500$ DERs are controlled; the limits for the active power of each device are  $[-50, 50]$ kW for energy storage resources and $[0, 50]$ kW for solar inverters. We consider the case where $p_{0,t}^{\mathrm{ref}}$ follows the trajectory shown in Figure~\ref{fig:regret_data}; real data with a granularity of one second is taken from~\cite{Dallanese2018feedback} to generate the non-controllable powers $\bw_t$, with the net power  $\ba_w^\top \bw_t$ plotted in Figure~\ref{fig:regret_data} as well. The sensitivity vector $\ba_w$ is computed as in~\cite{Dallanese2018feedback}. A Gaussian random variable with  zero mean and variance $10$ kW is utilized to generate the measurement error affecting $\hat{p}_{0,t}$.    

Figure~\ref{fig:regret_demand} illustrates the evolution of the regret $r_t$, averaged over $50$ experiments, in logarithmic scale. One can notice a linear decrease of the average regret during the first iteration of the algorithms; the regret then exhibit variations that are due to the considerable time-variability of the cost function (due to the large swings in the non-controllable powers $\bw_t$). The plot also provides a zoomed version (in linear scale), where the $3$-standard deviation confidence interval is also reported.

\section{Conclusions}
\label{sec:conclusions}
In this paper, we showed that cost function achieved by the online (proximal-)gradient method exhibits a linear convergence  to the optimal value functions within an error, for functions satisfying the  (proximal-)PL inequality, and when inexact gradient information is available. We derived bounds in expectation and in high probability, where for the latter we utilized a sub-Weibull model for the gradient errors. The convergence results are applicable to a number of learning and feedback-optimization tasks, where the cost functions may not be SC, but satisfies the PL inequality. Our results also provide new insights on the convergence of the (proximal-)gradient method for time-varying functions and with exact gradient information, and for the case of static optimization with inexact gradient information. The gradient error model is general, and it allows one to consider various sources of inaccuracy and gradient estimation techniques.

\bibliographystyle{IEEEtran}
\bibliography{biblio.bib}
\end{document}